\documentclass[11pt,reqno,UTF8,twoside]{amsart}
\usepackage{mathrsfs}
\usepackage{amsfonts,amssymb,amsmath,amsthm}
\usepackage[colorlinks=true,citecolor=red,linkcolor=blue]{hyperref}
\usepackage{titletoc}
\usepackage[noadjust]{cite}
\usepackage{geometry}
\usepackage{bm}
\usepackage{indentfirst}
\usepackage{graphicx}
\usepackage{stmaryrd}
\usepackage{tikz}
\usepackage{ulem}
\usepackage{color,soul}
\usepackage{setspace}
\usepackage[pagewise]{lineno} 
\usepackage{float}

\usepackage{todonotes}

\numberwithin{equation}{section}
\newtheorem{theorem}{Theorem}[section]
\newtheorem{lemma}[theorem]{Lemma}

\newtheorem{proposition}[theorem]{Proposition}
\newtheorem*{maintheorem}{Main Theorem}
\newtheorem{remark}[theorem]{Remark}
\theoremstyle{definition}
\newtheorem{example}[theorem]{Example}
\newtheorem{definition}[theorem]{Definition}

\newcommand{\Lip}{\operatorname{Lip}}

\newcommand{\dvg}{\operatorname{div}}
\newcommand{\Osc}{\operatorname{Osc}}
\newcommand{\supp}{\operatorname{supp}}
\newcommand{\good}{\operatorname{good}}
\newcommand{\bad}{\operatorname{bad}}

\makeatletter
\@namedef{subjclassname@2020}{\textup{2020} Mathematics Subject Classification}
\makeatother

\begin{document}
	\begin{abstract}
	    We study large deviations from the invariant measure for nonlinear Schr\"odinger equations with colored noises on determining modes. The proof is based on a new abstract criterion, inspired by [V. Jakši\'{c} et al., Comm. Pure Appl. Math., 68 (2015), 
        2108--2143]. To address the difficulty caused by fixed squeezing rate, we introduce a bootstrap argument to derive Lipschitz estimates for Feynman--Kac semigroups. This criterion is also applicable to wave equations and Navier--Stokes system. 
	\end{abstract}
    
    \title[Large deviation for Schr\"odinger equations]
	{\Large D\MakeLowercase{onsker--}V\MakeLowercase{aradhan large deviation principle for\\
    \vspace{2mm}locally damped and randomly forced} NLS \MakeLowercase{equations}}
	
\author[Y. Chen and S. Xiang]{Y\MakeLowercase{uxuan} C\MakeLowercase{hen}, S\MakeLowercase{hengquan} X\MakeLowercase{iang}}
		
\address[Yuxuan Chen]{School of Mathematical Sciences, Peking University, 100871, Beijing, China.}
\email{chen\underline{ }yuxuan@pku.edu.cn}

\address[Shengquan  Xiang]{School of Mathematical Sciences, Peking University, 100871, Beijing, China.}
\email{shengquan.xiang@math.pku.edu.cn}

	
	
	\maketitle
	
    
    \section{Introduction}

    The \textit{large deviation principle} (LDP) concerns the exponential decay of probabilities associated with rare events. It plays a crucial role in various fields, including statistical mechanics, dynamical systems and finance. For comprehensive overviews, see, e.g., the monographs \cite{Ell-85,DS89,DZ10}.
    
    Over the past two decades, there has been growing interest in the Donsker--Varadhan LDP for random evolutionary PDEs, with several research groups exploring this area \cite{Gou07-NS,Gou07-Burgers,JNPS-15,JNPS-18,JNPS-21,LL23,MN-18,Ner-19,NPX-23,CLXZ24}. Under appropriate assumptions on the noise structure, one can prove that the empirical distribution (or occupation measure) of the solution converges to the unique invariant measure. The Donsker--Varadhan LDP aims to describe its fluctuation at order one.

    \medskip
    
    While most existing studies have focused on parabolic PDEs, our motivation is to address two important facets that remain less understood: (1) hyperbolic equations; and (2) noises acting only on finitely many modes. These scenarios were investigated in our earlier joint work with Liu and Zhang \cite{CLXZ24}, which established the local LDP for wave equations. 
    
    The main contribution of the present work is to remove the local constraint for trajectories lying within the support of invariant measure. To achieve this, we establish a new general LDP criterion, inspired by \cite{JNPS-15}, which enables us to obtain the LDP for both hyperbolic and parabolic PDEs. The crux arises from the fixed squeezing rate, which creates obstacles for Lipschitz estimates of Feynman--Kac semigroups. We overcome this difficulty by introducing a bootstrap argument, which constitutes the technical novelty of this paper.

    \subsection{Main result}

    We formulate our main results for the nonlinear Schr\"odinger (NLS) equation, which simultaneously exhibits features of a hyperbolic equation and degenerate noise. Two more examples, the wave equation and the Navier–Stokes system, are provided in Section~\ref{subsec_application}. 

    \medskip
    
    The randomly forced nonlinear Schr\"odinger equation on 1D torus $\mathbb{T}:=\mathbb{R}/2\pi\mathbb{Z}$ reads
    \begin{equation}\label{NLS}
        \left\{\begin{array}{ll}
            iu_t+u_{xx}+ia(x)u=|u|^{p-1}u+\eta(t,x),\\
            u(0,x)=u_0(x),
        \end{array}\right.
    \end{equation}
    where the order $p\ge 3$ of nonlinearity is an odd integer. The function $a(\cdot)$ denotes the damping coefficient, which is smooth, non-negative, and not identically zero. The term $\eta(t,x)$ represents the random noise, which is spatially localized, degenerate in Fourier modes, bounded in a suitable Sobolev norm, and statistically $1$-periodic\footnote{For $T>0$, the random process $\eta(t)|_{t\ge 0}$ is called statistically $T$-periodic if $\eta(nT+t)|_{t\in [0,T)}$ are i.i.d.~for $n\in \mathbb{N}_0$. We specialize to $T=1$ for the sake of simplicity. All results in this paper remain true for arbitrary period $T>0$.}. See Section~\ref{subsec_NLS} for the precise noise structure. We emphasize that both the damping and the noise may vanish outside an open subset of $\mathbb{T}$.

    The system \eqref{NLS} is well-posed in the natural energy space $H^1(\mathbb{T})$. Since the noise is statistically $1$-periodic, the discrete-time sequence $u_n:=u(n)\, (n\in \mathbb{N}_0)$ forms a Markov process. Recently, the authors together with Zhang and Zhao \cite{CXZZ-25} proved that, under the assumptions of the Main Theorem below, this Markov process $(u_n,\mathbb{P}_u)$ admits a unique invariant measure $\mu_*\in \mathcal{P}(H^1)$; moreover, $\mu_*$ has compact support, and is exponentially mixing.

    A direct consequence of exponential mixing is the law of large numbers (see, e.g., \cite{LWXZZ24}): for any random initial distribution $u_0$, which we always assume to be independent of $\{\eta_n\}$, the empirical distribution of the process $u_n$ converges to $\mu_*$ almost surely, namely
    \[\frac{1}{n}(\delta_{u_1}+\cdots +\delta_{u_n})\rightharpoonup \mu_*,\quad a.s.\]
    Here $\delta_u$ stands for the Dirac measure concentrated on $u$.

    \medskip
    
    The Main Theorem provides the corresponding large deviation estimates. Recall that in the LDP theory, a good rate function refers to a lower semicontinuous function $I$ taking values in $[0,\infty]$, such that $\{I\le a\}$ is compact for any $a\in [0,\infty)$. The convention $\inf \emptyset=\infty$ is adopted.

    \begin{maintheorem}
        Let $\mathbf{B},\sigma>0$ be arbitrarily given. Assume the damping $a\colon \mathbb{T}\to [0,\infty)$ is smooth and not identically zero, and the noise $\eta$ can be formulated as in \eqref{noise-NLS}. Then there exists a constant $N\in \mathbb{N}$ such that if the non-negative numbers $b_{j,k}$ in \eqref{noise-NLS} satisfy 
        \begin{equation}\label{bounded-noise-0}
            \sum_{j\in\mathbb{N},\, k\in\mathbb{Z}} |b_{j,k}|^2 \langle k\rangle^{2(1+\sigma)}\leq \mathbf{B},
        \end{equation}
        and
        \begin{equation}\label{non-degenerate}
            b_{j,k}\neq 0\quad \text{for any }j,|k|\leq N,
        \end{equation}
        then there exists a convex good rate function $I\colon \mathcal{P}(\supp(\mu_*))\to [0,\infty]$, such that the LDP holds for the empirical distribution of NLS equation \eqref{NLS}, uniformly with respect to the initial distribution $\mathscr{D}(u_0)=\mu\in \mathcal{P}(\supp(\mu_*))$, in the following sense: for any Borel set $B\subset \mathcal{P}(\supp(\mu_*))$, 
        \begin{align*}
            -\inf_{\sigma\in B^\circ} I(\sigma) &\le \liminf_{n\to\infty}\inf_{\mu\in \mathcal{P}(\supp(\mu_*))}\frac{1}{n}\log \mathbb{P}_\mu \left(\frac{1}{n}\sum_{k=1}^n \delta_{u_k}\in B\right)\\&\le \limsup_{n\to \infty} \sup_{\mu\in \mathcal{P}(\supp(\mu_*))} \frac{1}{n}\log \mathbb{P}_\mu \left(\frac{1}{n}\sum_{k=1}^n \delta_{u_k}\in B\right)\le -\inf_{\sigma\in \overline{B}} I(\sigma),
        \end{align*}
        where $B^\circ$ and $\overline{B}$ refer to the interior and closure of $B$, respectively.
    \end{maintheorem}

    The proof is carried out in Section~\ref{subsec_NLS}, based on a general criterion formulated in Theorem~\ref{Thm LDP}. See Remark~\ref{Idef} for the definition of the rate function $I$. We also mention that the Main Theorem implies the LDP for test functions (or observables) via the contraction principle; see, e.g., \cite[Section 1.2]{JNPS-15}. The derivation of this corollary is standard and hence omitted.
    
    \subsection{Literature review}
    
    The study of Donsker--Varadhan LDP for stochastic parabolic PDEs goes back to Gourcy \cite{Gou07-NS,Gou07-Burgers}, based on earlier results for stochastic damped Hamiltonian system by Wu \cite{Wu01}. In these two papers, the LDP for the 1D viscous Burgers equation and 2D Navier--Stokes system was established, where the noise is white in time and sufficiently rough in space. 
    
    Then a general method for parabolic equations with non-degenerate (i.e.~acting on all modes) bounded kick force was introduced by Jakši\'{c}, Nersesyan, Pillet and Shirikyan \cite{JNPS-15}, using Kifer's criterion \cite{Kifer-90} and asymptotics of Feynman--Kac semigroups. This strategy was later extended to unbounded kicks \cite{JNPS-18}, white-in-time noise \cite{Ner-19}, and Lagrangian flows \cite{JNPS-21}.
    
    It is also worth mentioning that, for Navier--Stokes system on $\mathbb{T}^2$ with highly degenerate white-in-time noise, Malliavin calculus is a powerful tool: Hairer and Mattingly \cite{HM-06,HM-08} established mixing results, and recently Nersesyan, Peng and Xu \cite{NPX-23} proved the LDP. More recently, Liu and Lu \cite{LL23} further extended the LDP result to include a periodic external force.

    \medskip

    As for hyperbolic equations, due to the lack of smoothing effect, the LDP remains largely open. To the best of our knowledge, the only existing papers in this direction are \cite{MN-18} by Martirosyan and Nersesyan, and our earlier work with Liu and Zhang \cite{CLXZ24}, which proved local LDP for nonlinear wave equations with white-in-time and colored noises, respectively. Here the term ``local" means that the LDP lower bound only holds on a subset.

    \medskip

    Finally, we also refer the reader to \cite{CR04,SS06,CM10,BCF15,Mar17,FG-23,GGKS-23,OGB-24} for other LDP results for PDEs, and to \cite{Bourgain-96,BT-08,KMV-19,GOT-22,BDNY-24,GKO-24} for other related topics on random dispersive equations.

    \subsection{New ingredient}

    To prove the Main Theorem, we introduce a general criterion, modified from \cite{JNPS-15}, for the LDP in \textit{random dynamical systems} (RDS) of the form
    \[x_{n+1}=S(x_n,\eta_n),\quad x_0=x\in X.\]
    Here $X$ is a compact metric space, $\eta_n$ are i.i.d.~random variables, and the map $S$ dictates the evolution. The LDP criterion relies on two assumptions, which can be roughly described as:
    \begin{itemize}
        \item[\tiny$\bullet$] {\bf Irreducibility.} {\it The system can evolve from any initial state to any other.}
         
        \item[\tiny$\bullet$] {\bf Squeezing.} {\it Different trajectories may become closer at a fixed rate $q\in [0,1)$.}
    \end{itemize}
    See Section~\ref{subsec_RDS} for details. The key difficulty lies in the second assumption: the squeezing rate $q$ is fixed here, whereas in \cite{JNPS-15} it was allowed to be arbitrarily small. We show how this fixed rate complicates the Lipschitz estimate for Feynman–Kac semigroups, and explain our bootstrap argument that overcomes this issue.

    \medskip

    When demonstrating this criterion, the final task is to show (see \cite[Section 3.2]{JNPS-15})
    \begin{equation}\label{uf-1}
        |Q_n^V \mathbf{1}(x)-Q_n^V \mathbf{1}(x')|\lesssim \|Q_n^V \mathbf{1}\|_{\infty}d(x,x'),\quad \text{for any }x,x'\in X, n\in \mathbb{N},
    \end{equation}
    where $Q_n^V \mathbf{1}(x):=\mathbb{E}_x e^{V(x_1)+\cdots +V(x_n)}$ is the Feynman--Kac semigroup for $V\in L(X)$.
    To derive this Lipschitz estimate, \cite{JNPS-15} realized that the ``bad" events are where two trajectories stop squeezing at time $k<n$, whose contribution to the left-hand side of \eqref{uf-1} can be bounded by
    \[q^k e^{k\Osc(V)}\|Q_n^V \mathbf{1}\|_{L^\infty} d(x,x'),\]
    where $\Osc(V):=\sup V-\inf V$ refers to the oscillation of $V$. The summation over $k$ can be bounded only if $q$ is small enough so that $qe^{\Osc(V)}<1$. 

    \medskip
    
    Meanwhile, let us provide two examples to illustrate why the squeezing rate $q$ cannot be arbitrarily small in the cases of degenerate noise and hyperbolic equations. In both examples, the goal is to drive the solution to $0$ by choosing an appropriate control function $h$.

    \begin{example}[Parabolic equation with degenerate noise]
        Consider the linear heat equation 
        \[u_t-\Delta u=h\]
        on the torus $\mathbb{T}$. If $h$ is degenerate, namely $h$ contains only the Fourier modes $\{e^{ikx}:|k|<N\}$, then the high frequencies of the solution evolves as the same as when $h=0$. More precisely, for any $n\ge N$, we have
        \[\hat{u}(t,n)=e^{-n^2 t} \hat{u}_0(n).\]
        Thus the decay rate of the solution has a lower bound $e^{-N^2 t}$.

        We also mention that since $e^{-N^2 t}$ decreases to $0$ as $N$ tends to infinity, the squeezing rate can be arbitrarily small if $h$ is non-degenerate, by controlling sufficiently many modes using $h$ and letting the sufficiently high frequencies decay by linear evolution.
    \end{example}

    \begin{example}[Hyperbolic equation]
        Consider the damped Schr\"odinger equation 
        \[iu_t+\Delta u+iau=h\]
        on the torus $\mathbb{T}$, where $a>0$ is a constant damping. Unlike the heat equation, we have
        \[e^{(i\Delta -a)t}e^{inx}=e^{(-in^2-a)t}e^{inx},\]
        yielding the same decay rate $e^{-at}$ for all Fourier modes. If we control the first $N$ modes using $h$, and let the high frequencies decay by linear evolution, then the decay rate is always bounded from below by $e^{-at}$, independent of $N$.
    \end{example}

    Having clarified the obstacle of fixed squeezing rate, we now outline our new bootstrap argument for establishing Lipschitz estimates \eqref{uf-1}, which consists of two steps:

    \begin{itemize}
        \item[\tiny$\bullet$] We introduce a parameter $\lambda$ to characterize the exponential growth of $Q_n^V \mathbf{1}$, and work with the rescaled function $\lambda^{-n} Q_n^V \mathbf{1}$. By noticing that the conditional probabilities of bad events are also small, we can first establish a weaker Lipschitz estimate:
        \begin{equation}\label{uf-2}
            \lambda^{-n}|Q_n^V \mathbf{1}(x)-Q_n^V \mathbf{1}(x')|\lesssim \sup_{0\le k\le n} (\lambda^{-k}\|Q_k^V \mathbf{1}\|_\infty) d(x,x').
        \end{equation}

        \item[\tiny$\bullet$] The irreducibility implies that the dual semigroup $Q_n^{V*}$ admits an ``eigenmeasure" $\mu\in \mathcal{P}(X)$ with full support. If $\lambda^{-k} \|Q_k^V \mathbf{1}\|_\infty$ is unbounded, we can construct a non-negative function $F\not =0$ satisfying $\int_X F\, d\mu=0$, contradictory. We thus recover \eqref{uf-1} from \eqref{uf-2}.
    \end{itemize}

    It is noteworthy that irreducibility plays a role in Lipschitz estimates, while prior works only exploited squeezing property for this purpose. It is well-known that irreducibility is crucial to uniform LDP as well as the convexity of rate function; see Section~\ref{sec_discussion} for more discussions.

    \subsection{Organization}
    
    In Section~\ref{sec_framework} we establish the RDS criterion. Then we present three applications: nonlinear Schr\"odinger equation, nonlinear wave equation and Navier--Stokes system, all perturbed by degenerate colored noises. In Section~\ref{sec_ProofRDSLDP} we prove the criterion, based on the bootstrap argument explained above. Finally, Section~\ref{sec_discussion} contains some further discussions on the natural restriction $u_0\in \supp(\mu_*)$. In addition, the Appendix collects two auxiliary results.

    \subsection{Notations.} 
    
    $(\Omega,\mathcal F,\mathbb P)$ stands for a probability space, whose expectation is denoted by $\mathbb E$. 
    
    $(X,d)$ is a compact metric space. $\mathcal{B}(X)$ is the Borel $\sigma$-algebra. $B(x,r)$ is the open ball with center $x$ and radius $r$. The continuous functions constitute the Banach space $C(X)$, equipped with the supremum norm $\|\cdot\|_{\infty}$. The oscillation of $f$ is $\Osc(f):=\sup f-\inf f$. And $L(X)$ stands for Lipschitz functions, equipped with the Lipschitz norm $\|f\|_{L}:=\|f\|_{\infty}+ \Lip(f)$, where
    \[\Lip(f):= \sup\limits_{x\neq y}\frac{|f(x)-f(y)|}{d(x,y)}.\]
		
	$\mathcal{M}_+(X)$ and $\mathcal{P}(X)$ are the spaces of non-negative Borel measures and probability measures, both equipped with weak-* topology. $\mathcal{P}(X)$ is metrized by dual-Lipschitz distance:
    \[\|\mu-\nu\|_L^*:=\sup_{f\in L(X),\ \|f\|_L=1} |\langle f,\mu\rangle-\langle f,\nu\rangle|\quad \text{for }\mu,\nu\in \mathcal{P}(X).\]
    For $f\in C(X)$ and $\mu\in\mathcal{P}(X)$, we write $\langle f,\mu\rangle:=\int_X f\, d\mu$, and $\supp(\mu)$ for the support of $\mu$.

    \section{Model and results}\label{sec_framework}

    We first introduce an RDS model and state a general LDP criterion, which is a modification of \cite{JNPS-15}. As applications, we investigate three PDE examples, indicating that this criterion is suitable for a wide class of dissipative PDEs with degenerate bounded noise.
    
    \subsection{RDS model}\label{subsec_RDS}

    Consider the RDS on a compact metric space $(X,d)$, defined by
    \begin{equation}
        x_n=S(x_{n-1},\eta_{n-1}),\quad x_0=x\in X,\label{RDS}
    \end{equation}
    where $\eta_n$ are i.i.d.~$E$-valued random variables, and $E$ is a compact metric space. Furthermore, assume the map $S\colon X\times E\to X$ is continuous and uniformly Lipschitz in its first argument; that is, there exists a constant $C>0$, such that
    \[d(S(x,\zeta),S(x',\zeta))\le Cd(x,x'),\quad \text{for any }x,x'\in X,\ \zeta\in E.\]

    By virtue of the continuity of $S$ and the i.i.d.~nature of $\eta_n$, the RDS \eqref{RDS} defines a Feller family of discrete-time Markov processes $(x_n,\mathbb{P}_x)$; see, e.g., \cite[Section 1.3]{KS-12}. We denote the corresponding expected values by $\mathbb{E}_x$, and the Markov transition functions by $P_n(x,\cdot )$, i.e.,
    \[P_n(x,B):=\mathbb{P}_x(x_n\in B)\quad \text{for }B\in\mathcal{B}(X).\]
    The associated Markov semigroup $P_n\colon C(X)\rightarrow C(X)$ is given by
    \begin{equation*}
        P_nf(x):= \int_X f(y)\, P_n(x,dy).
    \end{equation*}
    And $P^*_n\colon \mathcal{P}(X)\rightarrow \mathcal{P}(X)$ refers to the dual semigroup. For the sake of simplicity, we write $P=P_1$ and $P^*=P_1^*$. A probability measure $\mu$ is called \textit{invariant} if $P^*\mu=\mu$.

    \newpage
    
    We impose the following two hypotheses for the general result:

    \begin{itemize}
        \item [$(\mathbf{I})$] {\bf Uniform irreducibility.} For any $\varepsilon>0$, there exists $N\in \mathbb{N}$ and $p>0$, such that
        \begin{equation}\label{irreducibility}
            P_N(x,B(y,\varepsilon))>p,\quad \text{for any } x,y\in X.
        \end{equation}
         
        \item [$(\mathbf{C})$] {\bf Coupling condition for squeezing.} There exists $q\in [0,1)$ and $C>0$ such that for any $x,x'\in X$, the pair $(P(x,\cdot),P(x',\cdot))$ admits a coupling\footnote{A coupling between $\mu,\nu\in\mathcal{P}(X)$ is an $X\times X$-valued random variable with marginal laws equal to $\mu$ and $\nu$.} $(R(x,x'),R'(x,x'))$ on a common probability space $(\Omega,\mathcal{F},\mathbb{P})$ satisfying
        \begin{equation}\label{squeezing}
            \mathbb{P}(d(R(x,x'),R'(x,x'))>q d(x,x'))\le Cd(x,x').
        \end{equation}
        In addition, the maps $R,R'\colon X\times X\times \Omega\rightarrow X$ are measurable.
    \end{itemize}
    
    It can be proved that $(\mathbf{I})$ and $(\mathbf{C})$ imply exponential mixing. Indeed, the following result is a special case of \cite[Theorem 3.1.7]{KS-12}.

    \begin{proposition}
        Under the above settings, suppose that hypotheses $(\mathbf{I})$ and $(\mathbf{C})$ are satisfied. Then the Markov process $(x_n,\mathbb{P}_x)$ defined by \eqref{RDS} admits a unique invariant measure $\mu_*\in\mathcal{P}(X)$. Moreover, $\supp(\mu_*)=X$, and there exist constants $C,\gamma>0$ such that
        \[\|P_n(x,\cdot)-\mu_*\|_L^*\le Ce^{-\gamma n},\quad \text{for any }x\in X,\ n\in \mathbb{N}.\]
    \end{proposition}

    For any random variable $x_0$ with law $\mu:=\mathscr{D}(x_0)\in \mathcal{P}(X)$, which we always assume to be independent of $\{\eta_n\}$, we introduce the empirical distributions
    \begin{equation}
        L_{n,\mu}:=\frac{1}{n}\sum_{k=1}^n\delta_{x_k}.\label{Lnxdef}
    \end{equation}
    
    \medskip

    Now we state the general LDP criterion, to be proved in the next section.
    
    \begin{theorem}\label{Thm LDP}
        Under the above settings, suppose that the hypotheses $(\mathbf{I})$ and $(\mathbf{C})$ are satisfied. Then there exists a convex good rate function $I\colon \mathcal{P}(X)\to [0,\infty]$ such that uniform LDP holds for empirical distributions. More precisely, for any Borel set $B\subset \mathcal{P}(X)$,
        \begin{align*}
            -\inf_{\sigma\in B^\circ} I(\sigma)&\le \liminf_{n\rightarrow\infty}\inf_{\mu\in \mathcal{P}(X)}\frac{1}{n}\log \mathbb{P}(L_{n,\mu}\in B)\\
            &\le \limsup_{n\rightarrow\infty}\sup_{\mu\in \mathcal{P}(X)} \frac{1}{n}\log \mathbb{P}(L_{n,\mu}\in B)\le -\inf_{\sigma\in \overline{B}} I(\sigma).
        \end{align*}
    \end{theorem}

    \begin{remark}\label{Idef}
        In Section~{\rm\ref{subsec_reduction}}, we characterize the rate function as a Legendre transform:
        \[I(\sigma):=\sup_{V\in C(X)} (\langle V,\sigma\rangle-\Lambda(V)),\]
        where $\Lambda\colon C(X)\to \mathbb{R}$ is the pressure function \eqref{def_Lambda}. Alternatively, the Donsker--Varadhan entropy formula {\rm\cite{DV-75}} holds (though we do not need this formulation in the present paper):
        \[I(\sigma)=-\inf_{f\in C(X),\ f>0} \int_X \log \frac{Pf}{f}\, d\sigma.\]
        (See, e.g., {\rm\cite[Section 3.2.2]{CLXZ24}} for the proof in a slightly more general setting.)
    \end{remark}

    \subsection{Proof of the Main Theorem}\label{subsec_NLS}

    As a first application of Theorem~\ref{Thm LDP}, in this subsection we prove the Main Theorem concerning the LDP for NLS equation \eqref{NLS}.

    \medskip
    
    To start with, let us state the precise setting on the noise structure. Denote by $\{\alpha_j:j\in\mathbb{N}\}$ an orthonormal basis of $L^2(0,1)$. We also define the trigonometric  basis of $L^2(\mathbb{T})$ by
    \[e_k(x)=\frac{1}{\sqrt{2\pi}} e^{ikx},\quad k\in \mathbb{Z}.\]
    Then $\{\alpha_j(t) e_k(x):j\in \mathbb{N},\ k\in \mathbb{Z}\}$ serves as an orthonormal basis of $L^2([0,1]\times \mathbb{T})$.

    \begin{itemize}
        \item ({\bf Noise structure}) The random noise $\eta(t,x)$ is specified as
        \begin{equation}\label{noise-NLS}
        \begin{aligned}
            &\eta(t,x)= \eta_n(t-n,x)\quad \text{for }t\in [n,n+1),\ n\in \mathbb{N}_0,\\
            &\displaystyle\eta_n(t,x)=\chi(x)\sum_{j\in\mathbb{N},\, k\in\mathbb{Z}}b_{j,k}(\theta^{n}_{j,k,1}+i\theta_{j,k,2}^{n})\alpha_j(t)e_k(x),\quad t\in[0,1).
        \end{aligned}
        \end{equation}
        Here $\chi\in C^\infty(\mathbb{T})\setminus \{0\}$ is an arbitrary non-zero function\footnote{\cite{CXZZ-25} assumed $\chi$ to be positive in an interval. But the proof still works out as long as $\chi$ is not identically zero.}, $b_{j,k}\ge 0$ are deterministic constants tending to $0$ sufficiently fast (see \eqref{bounded-noise-0}), and $\theta_{j,k,l}^{n}$ are independent real-valued random variables. Moreover, $\theta_{j,k,l}^{n}$ admits a probability density function $\rho_{j,k,l}$ supported by the interval $[-1,1]$, which is $C^1$ and satisfies $\rho_{j,k,l}(0)>0$.
    \end{itemize}
    
    In view of the boundedness condition \eqref{bounded-noise-0}, $\eta_n$ are i.i.d.~$L^2(0,1;H^{1+\sigma}(\mathbb{T}))$-valued random variables. Moreover, their common law $\mathscr{D}(\eta)\in \mathcal{P}(L^2(0,1;H^{1+\sigma}(\mathbb{T})))$ has compact support $E:=\supp(\mathscr{D}(\eta_n))$, due to a diagonal argument (or Tikhonov's theorem). We emphasize that $E$ is equipped with the topology of $L^2(0,1;H^{1+\sigma})$, rather than the weaker one of $L^2(0,1;H^1)$.
    
    \medskip
    
    Under the above settings, the deterministic version of \eqref{NLS}, i.e.~replacing $\eta$ by deterministic force $\zeta\in L^2(0,1;H^1(\mathbb{T}))$, is globally well-posed in $H^1(\mathbb{T})$ (see, e.g., \cite[Proposition A.6]{CXZZ-25}). We shall denote the time-$1$ solution map by
    \[S\colon H^1\times L^2(0,1;H^1)\to H^1,\quad S(u_0,\zeta)=u(1).\]
    Then $S$ is locally Lipschitz. Therefore, the random system \eqref{NLS} gives rise to an RDS by setting
    \[u_n=S(u_{n-1},\eta_{n-1}),\quad u_0\in H^1(\mathbb{T}).\]
    We also employ the self-explanatory notation $u_n=S_n(u_0;\eta_0,\dots ,\eta_{n-1})$.
    
    Recently, it has been proved in \cite[Main Theorem and Proposition 5.2]{CXZZ-25} that, under the assumptions of the Main Theorem, the system admits a unique invariant measure $\mu_*\in \mathcal{P}(H^1)$ and is exponentially mixing; moreover, $\supp(\mu_*)$ is a bounded subset of $H^{1+\sigma}$, and in particular compact in $H^1$. In fact, $\supp(\mu_*)$ coincides with the attainable set $\mathcal{A}(\{0\})$, which we now define.

    \begin{definition}[Attainable set]
        For $Z\subset H^1(\mathbb{T})$, we recursively define $\mathcal{A}_n(Z)$ by
        \begin{equation*}
            \mathcal{A}_n(Z):=\{S(u,\zeta):u\in \mathcal{A}_{n-1}(Z),\zeta\in E\},\quad n\in\mathbb{N}
        \end{equation*}
        with $\mathcal{A}_0(Z)=Z$. And the attainable set $\mathcal{A}(Z)$ is given by $\mathcal{A}(Z)=\overline{\bigcup_{n\in \mathbb{N}}\mathcal{A}_n(Z)}$ ($H^1$-closure).
    \end{definition}
    
    \medskip

    Now we are in a position to prove the Main Theorem via the general criterion.

    \begin{proof}[\rm \bf Proof of the Main Theorem] 
        The nonlinear Schr\"odinger equation \eqref{NLS} fits into the RDS model in Section~\ref{subsec_RDS}, by setting (the topologies are also inherited)
        \[X=\supp(\mu_*)\subset H^1(\mathbb{T}),\quad E=\supp(\eta_n)\subset L^2(0,1;H^{1+\sigma}(\mathbb{T})),\]
        and let $S$ to be the restriction of time-$1$ solution map to $X\times E$. As explained above, both $X$ and $E$ are compact, and $S$ is Lipschitz-continuous on $X\times E$. Moreover, $X$ is bounded in $H^{1+\sigma}(\mathbb{T})$. Once we can verify the assumptions in Theorem~\ref{Thm LDP}, then the Main Theorem is an immediately consequence. To this end, we recall some results from \cite{CXZZ-25}.
        
        Hypothesis $(\mathbf{C})$ has been demonstrated in \cite[Section 5.3]{CXZZ-25} via control theory, provided we equip $X$ with an equivalent norm (see \cite[Lemma 4.10]{CXZZ-25}). We mention that the non-degeneracy condition \eqref{non-degenerate} and the $L^2(0,1;H^{1+\sigma})$-compactness of $E$ are crucial to this point. 
        
        Meanwhile, a weaker version of $(\mathbf{I})$ is also derived in this prior work by global stability. Indeed, \cite[Section 5.3]{CXZZ-25} shows that \eqref{irreducibility} holds when $y=0$; that is, for any $\varepsilon_0>0$, one can find $N_0\in \mathbb{N}$ and $p_0>0$ such that
        \[P_{N_0}(x,B(0,\varepsilon_0))\ge p_0,\quad \text{for any }x\in X.\]
        On the other hand, for any $\varepsilon_1>0$, we have
        \[X=\mathcal{A}(\{0\})\subset \bigcup_{n=0}^\infty B_X(\mathcal{A}_n(\{0\}),\varepsilon_1),\]
        and the sequence $\mathcal{A}_n(\{0\})$ is increasing (since $S(0,0)=0$ and $0\in E$). By compactness we see that $X\subset B_X(\mathcal{A}_{N_1}(\{0\}),\varepsilon_1)$ for some $N_1\in \mathbb{N}$. In other words, for any $y\in X$, we can find $\zeta_0,\dots, \zeta_{N_1-1}\in E$ such that
        \[\|y-S_{N_1}(0;\zeta_0,\dots ,\zeta_{N_1-1})\|_X <\varepsilon_1.\]
        By uniform continuity, there exists $\delta>0$ (independent of $y\in X$), such that if $z\in X$ and $\xi_0,\dots ,\xi_{N_1-1}$ satisfy  $\|z\|_X<\delta$ and $\|\xi_i-\zeta_i\|_E <\delta$, then
        \[\|y-S_{N_1}(z;\xi_0,\dots ,\xi_{N_1-1})\|_X <2\varepsilon_1.\]
        Thus, since $\eta_n$ are i.i.d., we find that for any $z\in B_X(0,\delta)$,
        \begin{align*}
            P_{N_1}(z,B_X(y,2\varepsilon_1))&\ge \mathbb{P}(\|\eta_i-\zeta_i\|_E <\delta,\ i=0,\dots , N_1-1)\\
            &=\prod_{i=0}^{N_1-1}\mathbb{P}(\eta_i\in B_E (\zeta_i,\delta))\ge \left(\inf_{\zeta\in E} \mathbb{P}(\eta_n\in B_E(\zeta,\delta))\right)^{N_1}=:p_1.
        \end{align*}
        Since $\supp(\mathscr{D}(\eta_n))=E\ni \zeta\mapsto \mathbb{P}(\eta_n\in B_E(\zeta,\delta))$ is lower-semicontinuous and pointwise positive, and $E$ is compact, we conclude that $p_1$ is a positive constant. Finally, for any $\varepsilon>0$, let us set $\varepsilon_1=\varepsilon/2$, which determines the constant $\delta$, and then set $\varepsilon_0=\delta$. Then for $N:=N_0+N_1$ and any $x,y\in X$, by Markov property, we have
        \[P_N(x,B_X(y,\varepsilon))\ge P_{N_0}(x,B_X(0,\delta))\inf_{z\in B_X(0,\delta)} P_{N_1}(z,B_X(y,\varepsilon))\ge p_0 p_1>0.\]
        Thus the hypothesis $(\mathbf{I})$ is valid. Now the proof is complete.
    \end{proof}
    
    \subsection{Two more applications}\label{subsec_application}
    
    We present two more applications to random PDEs. These examples indicate that our criterion is suitable for both hyperbolic and parabolic PDEs.

    \subsubsection{Nonlinear wave equation}\label{subsubsec_NLW}
    
    The wave equation in consideration is
     \begin{equation}
        \begin{cases}
            \boxempty u+a(x)\partial_t u+u^3=\eta(t,x),\quad x\in D,\\
            u|_{\partial D}=0,\quad u[0]=(u_0,u_1)=:\boldsymbol{u}.
        \end{cases}\label{waveequation}
    \end{equation} 
    Here $D$ is a smooth bounded domain in $\mathbb{R}^3$, with boundary $\partial D$ and outer normal vector $n(x)$, the notation $\boxempty:=\partial_{tt}^2-\Delta$ represents the d'Alembert operator, and $u[t]:=(u,\partial_t u)(t)$. In order to formulate our settings for damping $a(x)$ and random noise $\eta(t,x)$, let us first recall the notion of $\Gamma$-type domain, initially used by Lions \cite{Lions-88}.
    
    \begin{definition}[$\Gamma$-type domain]\label{Def-Gamma}
        A $\Gamma$-type domain is a subdomain of $D$ in the form
        \[N_\delta(x_0):=\{x\in D:|x-y|<\delta\text{ for some }y\in\Gamma(x_0)\},\]
        where $x_0\in\mathbb R^3\setminus \overline D$, $\delta>0$ and $\Gamma(x_0)=\{x\in\partial D:(x-x_0)\cdot n(x)>0\}$. 
    \end{definition}
    
    \noindent This geometric setting is involved both in the localization of $a(x)$ and the structure of $\eta(t,x)$:

    \begin{itemize}
        \item ({\bf Localized structure}) {\it The function $a(\cdot)\in C^\infty(\overline{D})$ is non-negative, and there exists a $\Gamma$-type domain $N_\delta(x_0)$ and a constant $a_0>0$ such that
        \[a(x)\geq a_0,\quad \text{for any }x\in N_\delta(x_0).\]
        Moreover, let $\chi(\cdot)\in C^\infty(\overline{D})$ satisfy that for a $\Gamma$-type domain $N_{\delta'}(x_1)$ and $\chi_0>0$,
        \[\chi(x)\geq \chi_0,\quad\text{for any }x\in  N_{\delta'}(x_1).\]}
        
        \item ({\bf Noise structure}) {\it The random noise $\eta(t,x)$ is specified as
        \begin{equation}\label{noise-NLW}
        \begin{aligned}
            &\eta(t,x)= \eta_n(t-n,x)\quad \text{for }t\in [n,n+1),\ n\in\mathbb{N}_0,\\ 
            &\displaystyle\eta_n(t,x)=\chi(x)\sum_{j,k\in\mathbb{N}}b_{jk}\theta^{n}_{jk}\alpha_j(t)e_k(x)\quad \text{for }t\in [0,1),
        \end{aligned}
        \end{equation}
        where $b_{jk}\ge 0$ are constants, $\theta_{jk}^{n}$ are independent real random variables such that $\theta_{jk}^{n}$ admits $C^1$ density $\rho_{jk}$, which is supported by $[-1,1]$ and satisfy $\rho_{jk}(0)>0$; the sequence $\{\alpha_j:j\in\mathbb{N}\}$ denotes a smooth orthonormal basis of $L^2(0,1)$; the sequence $\{e_k:k\in \mathbb{N}\}$ are the eigenfunctions of Dirichlet Laplacian $-\Delta$ on $D$, forming an orthonormal basis of $L^2(D)$. The eigenvalue corresponding to $e_k$ is denoted with $\lambda_k$.}
    \end{itemize}

    \vspace{3mm}

    The phase space for \eqref{waveequation} is the energy space $\mathcal H:=H_0^1(D)\times L^2(D)$. Under the above settings, the integer-time solutions $\boldsymbol{u}_n:=u[n]$ form a Markov process. The uniqueness of invariant measure $\mu_*$ and exponential mixing has recently been established in \cite{LWXZZ24}.\footnote{Strictly speaking, \cite{LWXZZ24} only proved exponential mixing for $T$-periodic noises when $T$ is sufficiently large. In fact, the largeness of $T$ is only used to ensure $\|S_a(T)\|_{\mathcal{H}\to \mathcal{H}}<1$, where $S_a(t)$ is the operator semigroup generated by $\Delta-a(x)\partial_t$. This can be overcome, for any $T>0$, by selecting an equivalent norm on $\mathcal{H}$; see \cite[footnote 2]{CXZZ-25}.} As in the previous example, this Markov process fits into the RDS model \eqref{RDS}, by setting $S$ as the time-$1$ solution map and restrict $S$ to $\supp(\mu_*)\times \supp(\mathscr{D}(\eta_n))$. Another application of Theorem~\ref{Thm LDP} yields:
    
    \begin{theorem}\label{thm_wave}
        Let $\mathbf{B}>0$ be arbitrarily given. Assume the above settings are valid. Then there exists a constant $N\in \mathbb{N}$, such that if the constants $b_{jk}$ in \eqref{noise-NLW} satisfy
        \begin{equation}
            \sum_{j,k\in\mathbb{N}}b_{jk}\lambda_k^{2/7}\|\alpha_j\|_{_{L^\infty(0,1)}} \leq \mathbf{B}\quad\text{and}\quad  b_{jk}>0\ \text{ for }\ 1\leq j,k\leq N,\label{noise-structure2}
        \end{equation}
        then for any initial distribution $\mu\in \mathcal{P}(\supp(\mu_*))$, the empirical distributions satisfy uniform LDP with convex good rate function, in the sense of Theorem~{\rm\ref{Thm LDP}}.
    \end{theorem}

    \begin{proof}[\rm \bf Proof of Theorem~\ref{thm_wave}]
        As in the proof of the Main Theorem, the hypothesis $(\mathbf{C})$ and a weaker version of $(\mathbf{I})$ has been established in \cite{LWXZZ24}. In order to recover the actual hypothesis $(\mathbf{I})$, we repeat the previous arguments for the Main Theorem, as they only depend on the RDS structure and Markov property, without relying on specific PDE settings. Now the proof is complete.
    \end{proof}
   
    \subsubsection{Navier--Stokes system}\label{subsubsec_NS}

    Finally, we investigate the 2D Navier--Stokes system with bounded kick force. In particular, our result extends \cite[Theorem A]{JNPS-15} by allowing the noise to degenerate in high frequencies. The system in consideration is studied in, e.g., \cite{KS-12,JNPS-15,JNPS-18}, which reads
    \begin{equation}
        \begin{cases}
            \partial_t u-\nu \Delta u+u\cdot \nabla u+\nabla p=\eta(t,x),\quad\dvg u=0,\quad x\in D,\\
            u|_{\partial D}=0,\quad u(0)=u_0.
        \end{cases}\label{NSsystem}
    \end{equation}
    Here $D$ is a smooth bounded domain in $\mathbb{R}^2$, with boundary $\partial D$ and outer normal vector $n(x)$, the constant $\nu>0$ is the viscosity, $u=(u^1,u^2)$ is the velocity field, and $p$ is the pressure. The phase space for \eqref{NSsystem} is specified as 
    \[H:=\{u\in L^2(D;\mathbb{R}^2): \dvg u=0\text{ in }D,\ u\cdot n=0\text{ on }\partial D\}.\]
    The noise structure can be formulated as follows:

    \begin{itemize}
        \item ({\bf Noise structure}) {\it The random noise $\eta$ is a bounded kick of the form
        \begin{equation}\label{noise_NS}
            \eta(t,x)=\sum_{n=1}^\infty \delta(t-n)\eta_{n-1}(x),\quad \eta_n(x)=\sum_{k=1}^\infty b_k \theta_k^n \varphi_k(x),
        \end{equation}
        where $\delta$ is the Dirac measure, $b_k\ge 0$ are constants, $\theta_k^n$ are independent real random variables with $C^1$ density $\rho_k$, which is supported by $[-1,1]$ and satisfies $\rho_k(0)>0$; and $\{\varphi_k:k\in \mathbb{N}\}$ is an orthonormal basis of $H$.
    }
    \end{itemize}
    
    The solution of \eqref{NSsystem} is defined to be right-continuous at time $t=n$. Then $u_n:=u(n)$ gives rise to a Markov process. Equivalently, let $S_1\colon H\to H$ be the time-$1$ solution map for the unforced system (i.e.~\eqref{NSsystem} with $\eta=0$), and define $S(u,\eta)=S_1(u)+\eta$. Then we have
    \[u_n=S(u_{n-1},\eta_{n-1}),\quad u_0\in H.\]
    
    The uniqueness of invariant measure $\mu_*$ and exponential mixing for this system can be found in, e.g., \cite[Theorem 3.2.9]{KS-12}. Now we state the corresponding LDP result.
    
     \begin{theorem}\label{thm_2Dns}

        Let $B>0$ be arbitrarily given. Assume the above settings are valid. Then there exists a constant $N\in \mathbb{N}$ such that if the constants $b_k$ in \eqref{noise_NS} satisfy
        \[\sum_{k=1}^\infty b_k^2\le B_0\quad \text{and}\quad b_k>0\ \text{ for }\ 1\le k\le N,\]
        then for any initial distribution $\mu\in \mathcal{P}(\supp(\mu_*))$, the empirical distributions satisfy uniform LDP with convex good rate function, in the sense of Theorem~{\rm\ref{Thm LDP}}.
    \end{theorem}

    \begin{proof}[\rm \bf Proof of Theorem~\ref{thm_2Dns}]
        Once again, we restrict $S$ to $\supp(\mu_*)\times \supp(\mathscr{D}(\eta_n))$, which suits into the RDS model in Section~\ref{subsec_RDS}. For Navier--Stokes system, both hypotheses $(\mathbf{I})$ and $(\mathbf{C})$ can be found in the literature (see, e.g., \cite[Section 1.3]{JNPS-15}\footnote{Indeed, \cite{JNPS-15} invoked a stronger coupling condition, which implies $(\mathbf{C})$ in general; see, e.g., \cite[Section 3.2.2]{KS-12}.}). Therefore, the proof is automatic.
    \end{proof}

    Finally, we mention a similar application to Navier--Stokes system perturbed by degenerate space-time localized bounded noise, which has been studied in \cite{Shi-15} for exponential mixing. In this case the noise is akin to our setting for Schr\"odinger and wave equations. For the sake of simplicity, we omit the detailed description of this example, and claim that the LDP holds for the empirical distributions, as both $(\mathbf{I})$ and $(\mathbf{C}$) has essentially been verified in that work.

    \section{Proof for the LDP criterion}\label{sec_ProofRDSLDP}
    
    This section is devoted to demonstrating the abstract LDP criterion formulated in Theorem~\ref{Thm LDP}. We first reduce the issue to Lipschitz estimates of Feynman--Kac semigroups (see Proposition~\ref{prop_UF}). Then we finish this task via a new bootstrap argument. Throughout this section, we work in the general RDS setting in Section~\ref{subsec_RDS}. 

    \subsection{Some reductions}\label{subsec_reduction}

    We carry out three stages of reductions, following the route of \cite{JNPS-15}: (1)~Kifer's LDP criterion \cite{Kifer-90} for random probability measures; (2)~asymptotics of Feynman--Kac semigroups \cite{JNPS-15}; and (3)~Lipschitz estimates. No originality is claimed through these reductions.

    \vspace{3mm}

    \textit{Reduction 1. Kifer's criterion.} In order to investigate the uniform LDP, let us consider the directed set\footnote{A directed set $(\Theta,\prec)$ is a non-empty set $\Theta$ equipped with a pre-order $\prec$ (i.e.~a reflexive and transitive binary relation), such that any finite subset admits an upper bound. A net $\{x_\theta:\theta\in \Theta\}$ in a topological space $X$ refers to elements labeled by this directed set. The net $x_\theta$ converges to $x$, if for any neighborhood $U$ of $x$, there exists $\theta_0\in \Theta$ such that $x_\theta\in U$ whenever $\theta_0\prec \theta$.} $\Theta:=\mathbb{N}\times \mathcal{P}(X)$, endowed with the pre-order $\prec$ defined by 
    \begin{equation*}
        (n_1,\mu_1)\prec (n_2,\mu_2)\quad\text{if and only if}\quad n_1\le n_2.
    \end{equation*} 
    These notations are compatible with Appendix~\hyperlink{Appendix_A}{A}. Then Theorem~\ref{Thm LDP} would be a consequence of Kifer's criterion (Theorem~\ref{Thm Kifer}), once the assumptions therein are valid.
    
    To illustrate, set $r_\theta=n$ for $\theta=(n,\mu)\in \Theta$, and recall that the empirical distributions $L_\theta=L_{n,\mu}$ are defined by \eqref{Lnxdef}. Assume for now that, for any $V\in C(X)$, the right-hand side of
    \begin{equation}\label{def_Lambda}
        \Lambda(V):=\lim_{\theta\in \Theta} \frac{1}{r_\theta} \log \mathbb{E} e^{r_\theta \langle V,L_\theta\rangle}
    \end{equation}
    converges (alternatively, use $\limsup$ instead of $\lim$ to avoid the issue of convergence). We point out that $\Lambda(V)$ does not depend on the initial distribution $\mu$. More precisely, owing to
    \[\frac{1}{r_\theta} \log \mathbb{E} e^{r_\theta \langle V,L_\theta\rangle}=\frac{1}{n}\log \int_X  \mathbb{E}_x e^{V(x_1)+\cdots +V(x_n)}\, \mu(dx),\]
    and $\delta_x\in \mathcal{P}(X)\, (\forall x\in X)$, 
    the assumption above is equivalent to that, the functions
    \[X\ni x\mapsto \frac{1}{n} \log \mathbb{E}_x e^{V(x_1)+\cdots +V(x_n)}\]
    converge uniformly on $X$ to a constant $\Lambda(V)$. In the LDP theory, the map $\Lambda\colon C(X)\to \mathbb{R}$ is called the pressure function, which is readily seen to be $1$-Lipschitz continuous and convex.
    
    We now recall some standard facts from convex analysis. The Legendre transform of $\Lambda$ refers to the map $I\colon \mathcal{P}(X)\to [0,\infty]$ defined by
    \begin{equation}
        I(\sigma):=\sup_{V\in C(X)} (\langle V,\sigma\rangle-\Lambda(V))\quad \text{for }\sigma\in \mathcal{P}(X).\label{legendre}
    \end{equation}
    which is lower semicontinuous and convex. Moreover, $\Lambda$ can be reconstructed by the dual relation:
    \begin{equation}
        \Lambda(V)=\sup_{\sigma\in \mathcal{P}(X)} (\langle V,\sigma\rangle-I(\sigma)).\label{dualrelation}
    \end{equation}
    Note that $I$ is a good rate function, according to lower semicontinuity and compactness of $X$. 
    
    \begin{definition}[Equilibrium state]\label{def_es}
        Let $\Lambda\colon C(X)\to \mathbb{R}$ and $I\colon \mathcal{P}(X)\to [0,\infty]$ be defined as above. Given $V\in C(X)$, if the supremum on the right-hand side of \eqref{dualrelation} is attained at $\sigma_V\in \mathcal{P}(X)$, then $\sigma_V$ is called an equilibrium state for $V$.
    \end{definition}

    In order to apply Theorem~\ref{Thm Kifer}, it suffices to demonstrate the following lemma:
    
    \begin{lemma}\label{lemma_verifyKifer}
        Assume the hypotheses $(\mathbf{I})$ and $(\mathbf{C})$ are satisfied, then $\Lambda\colon C(X)\to \mathbb{R}$ is well-defined by \eqref{def_Lambda}, and the equilibrium state is unique for any $V\in L(X)$.
    \end{lemma}

    \begin{proof}[\rm \bf Proof of Theorem~\ref{Thm LDP}]
        Since $L(X)$ is a dense linear subspace of $C(X)$, Lemma~\ref{lemma_verifyKifer} implies that Theorem~\ref{Thm Kifer} is applicable, which immediately yields Theorem~\ref{Thm LDP}.
    \end{proof}

    \vspace{3mm}
    \textit{Reduction 2. Feynman--Kac semigroups.} The content of Lemma~\ref{lemma_verifyKifer} is closely related to the asymptotics of Feynman--Kac semigroups, which we now introduce.
    
    \begin{definition}[Feynman--Kac semigroup]
        For any $V\in C(X)$ and $n\in \mathbb{N}$, define
        \[Q^{V}_nf(x):=\mathbb{E}_x f(x_n)e^{V(x_1)+\cdots +V(x_n)}\quad\text{for }f\in C(X).\]
        It is easy to see $Q^{V}_n\colon C(X)\to C(X)$ forms a semigroup, owing to the Markov property. Denote its dual semigroup with $Q_n^{V*}\colon \mathcal{M}_+(X)\rightarrow\mathcal{M}_+(X)$.
    \end{definition}
    
    For simplicity, write $Q^V$ and $Q^{V*}$ in place of $Q_1^V$ and $Q^{V*}_1$.

    \begin{lemma}\label{lemma_FK}
        Assume the hypotheses $(\mathbf{I})$ and $(\mathbf{C})$ are satisfied, and $V\in L(X)$. Then there exists $\lambda_V>0$, $h_V\in C(X)$ and $\mu_V\in\mathcal{P}(X)$, such that for any $f\in C(X)$ and $\nu\in\mathcal{M}_+(X)$,
        \begin{align*}
            &Q^{V}h_V=\lambda_Vh_V,\quad  Q^{V*}\mu_V=\lambda_V\mu_V,\quad\langle h_V,\mu_V\rangle=1,\quad h_V>0,\quad \supp(\mu_V)=X,\\
            &\lambda_V^{-n}Q_n^{V}f\rightarrow\langle f,\mu_V\rangle h_V\quad\text{ in }C(X),\quad \lambda_V^{-n}Q_n^{V*}\nu\rightarrow\langle h_V,\nu\rangle \mu_V\quad\text{ in }\mathcal{M}_+(X).
        \end{align*}
    \end{lemma}

    \begin{proof}[\rm \bf Proof of Lemma~\ref{lemma_verifyKifer}]
        It follows from Lemma~\ref{lemma_FK} verbatim as in \cite[Section 4]{JNPS-15}.
    \end{proof}

    \vspace{3mm}
    \textit{Reduction 3. Lipschitz estimates.} The operators $Q_n^V$ form a generalized Markov semigroup, introduced in Appendix~\hyperlink{Appendix_B}{B}. Then Lemma~\ref{lemma_FK} is a consequence of Theorem~\ref{Thm Feynman-Kac}. Once again we need to verify the assumptions therein: uniform irreducibility and uniform Feller property. In particular, the latter relies on the following Lipschitz estimate.

    \begin{proposition}\label{prop_UF}
        Assume the hypotheses $(\mathbf{I})$ and $(\mathbf{C})$ are satisfied. Given $V\in L(X)$, there exists a constant $C_V<\infty$ such that for any $f\in L(X)$, $x,x'\in X$ and $n\in \mathbb{N}$,
        \begin{equation}
            |Q_n^V f(x)-Q_n^V f(x')|\le C_V\|f\|_L \|Q_n^V \mathbf{1}\|_\infty d(x,x').\label{UF inequality}
        \end{equation}
    \end{proposition}

    \begin{proof}[\rm \bf Proof of Lemma~\ref{lemma_FK}]
        It suffices to verify the two assumptions in Theorem~\ref{Thm Feynman-Kac}. The uniform irreducibility for $Q_n^V$ follows from $Q_n^V(x,\cdot)\ge e^{n\inf V} P_n(x,\cdot)$ and hypothesis $(\mathbf{I})$. The uniform Feller property is a consequence of Proposition~\ref{prop_UF} as follows. Set $\mathcal{C}:=L(X)\cap \{f\in C(X):f>0\}$. Then $\mathcal{C}$ is determining (see Definition~\ref{def_determine}), and $\|Q_n^V \mathbf{1}\|_\infty\le (\inf f)^{-1} \|Q_n^V f\|_\infty$ for any $f\in \mathcal{C}$. Thanks to \eqref{UF inequality}, we find that $\|Q_n^V f\|_\infty^{-1} Q_n^V f$ is equicontinuous. This completes the proof.
    \end{proof}

    \subsection{Lipschitz estimates for $Q_n^V$}\label{subsec_uf}
    
    It remains to address Proposition~\ref{prop_UF}. We need to estimate
    \begin{equation}
        Q_n^V f(x)-Q_n^V f(x').\label{QnVf-QnVf'}
    \end{equation}
    To this end, we bootstrap between two sequences: the rescaled Lipschitz constant
    \begin{equation}\label{def_Ln}
        L_n:=\lambda_1^{-n}\Lip(Q_n^V f),
    \end{equation}
    where $\lambda_1>0$ is a constant to be fixed later; and the rescaled $L^\infty$ norm
    \begin{equation}\label{def_Mn}
        M_n:=\sup_{0\le k\le n} \lambda_1^{-k}\|Q_k^V \mathbf{1}\|_\infty.
    \end{equation}
    We will simultaneously obtain the boundedness of $L_n$ and $M_n$ in two steps.

    \medskip
    
    \textit{Step 1: $L_n\le CM_n$.} We start with the scale $\lambda_1$ capturing the exponential growth of $\|Q_n^V \mathbf{1}\|_\infty$.
    
    \begin{lemma}\label{lemma_lambda1}
        For $V\in L(X)$, there exists a constant $\lambda_1>0$, such that for any $n\in \mathbb{N}$,
        \begin{equation}
        \|Q_n^V \mathbf{1}\|_\infty\ge \lambda_1^n.\label{Mn>=lambdan}
        \end{equation}
        Moreover, for any $\rho\in [0,1)$, the following positive series is convergent:
        \begin{equation}
        \sum_{n=0}^\infty \frac{\|Q_n^V\mathbf{1}\|_\infty}{\lambda_1^n} \rho^n<\infty.\label{geoserie}
        \end{equation}
    \end{lemma}

    \begin{proof}[\rm \bf Proof of Lemma~\ref{lemma_lambda1}]
        Since $Q_n^V$ is a semigroup and $\|Q_n^V f\|_\infty \le \|f\|_\infty \|Q_n^V \mathbf{1}\|_\infty$, we have
        \[\|Q_{m+n}^V\mathbf{1}\|_{\infty}\le \|Q_m^V\mathbf{1}\|_{\infty}\|Q_n^V\mathbf{1}\|_{\infty}.\]
        Owing to Fekete's Lemma on subadditive sequences, we find that $\frac{1}{n}\log \|Q_n^V \mathbf{1}\|_\infty$ converges to its infimum. As $e^{n\inf V}\le Q_n^V \mathbf{1}\le e^{n\sup V}$ by definition, there exists $\lambda_1\in [e^{\inf V},e^{\sup V}]$ such that
        \begin{equation}
            \lim_{n\to \infty} \frac{1}{n}\log \|Q_n^V \mathbf{1}\|_{\infty}=\log \lambda_1=\inf_{n\in \mathbb{N}} \frac{1}{n}\log \|Q_n^V \mathbf{1}\|_\infty.\label{lambda1definition}
        \end{equation}
        The second equality implies \eqref{Mn>=lambdan}, and the first equality implies \eqref{geoserie} by root test.
    \end{proof}

    In order to estimate \eqref{QnVf-QnVf'}, we exploit hypothesis $(\mathbf{C})$ to build up a coupling. Thanks to \eqref{squeezing}, for $x,x'\in X$ we can define a coupling operator on $X\times X$ by the relation
    \begin{equation*}
        \mathcal{R}(x,x'):=  (R(x,x'),R'(x,x')).
    \end{equation*}
    Let $\{(\Omega_n,\mathcal{F}_n,\mathbb{P}_n):n\in \mathbb{N}_0\}$ be a sequence of i.i.d.~copies of the probability space on which $\mathcal{R}$ is defined. Let $(\boldsymbol{\Omega},\boldsymbol{\mathcal{F}},\mathbb{P})$ be the product probability space of $\{(\Omega_n,\mathcal{F}_n,\mathbb{P}_n):n\in \mathbb{N}_0\}$. For any $x,x'\in X$ and $\boldsymbol{\omega}\in\boldsymbol{\Omega}$, we recursively define $\{(x_n,x_n'):n\in \mathbb{N}\}$ by 
    \begin{align*}
        (x_n(\boldsymbol{\omega}),x_n'(\boldsymbol{\omega}))&=\mathcal{R}(x_{n-1},x_{n-1}',\omega_{n-1}),
    \end{align*}
    where $(x_0,x_0'):=(x,x')$. In particular, by induction one finds that the laws of $x_n$ and $x_n'$ coincide with $P_n(x,\cdot)$ and  $P_n(x',\cdot)$, respectively. Let us denote by $\mathbb{P}_{\boldsymbol{x}}$ the law of the process issuing from $\boldsymbol{x}=(x,x')$, by $\mathbb{E}_{\boldsymbol{x}}$ the corresponding expectation, and by $\boldsymbol{\mathcal{F}}_n$ the natural filtration. 

    For any $k\in \mathbb{N}_0$, we introduce the following events (see \cite[Section 3]{JNPS-15})
    \begin{align*}
        A_k=\{d(x_{k+1},x_{k+1}')\le qd(x_k,x_k')\},\quad B_k=\bigcap_{j=0}^{k-1} A_j,\quad C_k=\bigcap_{j=0}^{k-1} A_j \cap A_k^C.        
    \end{align*}
    Then for any $n\in \mathbb{N}$, the sample space $\boldsymbol{\Omega}$ can be partitioned into 
    \begin{equation}
        \boldsymbol{\Omega}=\bigcup_{k=0}^{n-1} C_k \cup B_n.\label{Omegadecompose}
    \end{equation}
    The event $B_n$ is ``good", as we have $d(x_k,x_k')\le q^k d(x,x')$ for all $0\le k\le n$. Meanwhile, the ``bad" event $C_k$ is small, for example $\mathbb{P}(C_k)\le Cq^k d(x,x')$ as observed in \cite{JNPS-15}.
    Nevertheless, we need a finer estimate for its conditional probability.

    \begin{lemma}\label{lem_conditionprob}
        Assume the hypothesis $(\mathbf{C})$ is satisfied, then almost surely
        \begin{equation}
            \mathbb{E}(\mathbf{1}_{C_k}|\boldsymbol{\mathcal{F}}_k)\le Cq^k d(x,x'),\quad \text{for any }k\in \mathbb{N}_0.\label{conditionprob}
        \end{equation}
    \end{lemma}

    We emphasize that the constant $q\in [0,1)$ is the fixed squeezing rate from hypothesis $(\mathbf{C})$.

    \begin{proof}[\rm \bf Proof of Lemma~\ref{lem_conditionprob}]
        In view of $C_k=B_k\cap A_k^C$ and $B_k\in \boldsymbol{\mathcal{F}}_k$, we get
        \[\mathbb{E}(\mathbf{1}_{C_k}|\boldsymbol{\mathcal{F}}_k)=\mathbb{E}(\mathbf{1}_{B_k}\mathbf{1}_{A_k^C}|\boldsymbol{\mathcal{F}}_k)=\mathbf{1}_{B_k}\mathbb{E}(\mathbf{1}_{A_k^C}|\boldsymbol{\mathcal{F}}_k).\]
        Thanks to \eqref{squeezing} and our construction of the coupling, almost surely
        \begin{equation}
            \mathbb{E}(\mathbf{1}_{A_k^C}|\boldsymbol{\mathcal{F}}_k)=\mathbb{P}(d(x_{k+1},x_{k+1}')>q d(x_k,x_k')|\boldsymbol{\mathcal{F}}_k)\le Cd(x_k,x_k'),\label{squeezingattimen}
        \end{equation}
        And as $d(x_k,x_k')\le q^kd(x,x')$ on $B_k$, we conclude that
        \[\mathbb{E}(\mathbf{1}_{C_k}|\boldsymbol{\mathcal{F}}_k)\le C\mathbf{1}_{B_k} d(x_k,x_k')\le C\mathbf{1}_{B_k}q^k d(x,x')\le Cq^k d(x,x').\qedhere\]
    \end{proof}

    Now we are in a position to finish Step 1.

    \begin{lemma}\label{lem_Ln<=Mn}
        Assume the hypothesis $(\mathbf{C})$ is satisfied, then for any $V\in L(X)$, there exists a constant $C_V>0$, such that for any $f\in L(X)$, the rescaled Lipschitz constant $L_n$ and rescaled $L^\infty$ norm $M_n$ defined by \eqref{def_Ln} and \eqref{def_Mn} satisfy
        \begin{equation}
            L_n\le C_V\|f\|_L M_n,\quad \text{for any }n\in \mathbb{N}.\label{Ln<=Mninequality}
        \end{equation}
    \end{lemma}

    \begin{proof}[\rm \bf Proof of Lemma~\ref{lem_Ln<=Mn}]
        In view of the partition \eqref{Omegadecompose}, we decompose \eqref{QnVf-QnVf'} as
        \begin{equation}
            Q^V_n f(x)-Q^V_n f(x')=\sum_{k=0}^{n-1} J_n^{\bad,k}+J_n^{\good},\label{UF inequality1}
        \end{equation}
        where
        \begin{align*}
            &J_n^{\bad,k}:=\mathbb{E}_{\boldsymbol{x}}( \mathbf{1}_{C_k}(f(x_n)e^{\sum_{i=1}^n V(x_i)}-f(x_n')e^{\sum_{i=1}^n V(x_i')})),\\
            &J_n^{\good}:=\mathbb{E}_{\boldsymbol{x}}(\mathbf{1}_{B_n}(f(x_n)e^{\sum_{i=1}^n V(x_i)}-f(x_n')e^{\sum_{i=1}^n V(x_i')})).
        \end{align*}

    \medskip

    First we estimate $J_n^{\bad,k}$, exploiting the smallness of $C_k$. Using Markov property, we have 
    \begin{align*}
         \left|\mathbb{E}_{\boldsymbol{x}} (\mathbf{1}_{C_k}f(x_n)e^{\sum_{i=1}^n V(x_i)})\right|&\le \|f\|_\infty  \mathbb{E}_{\boldsymbol{x}} (\mathbf{1}_{C_k}e^{\sum_{i=1}^{k+1} V(x_i)}\mathbb{E}_{\boldsymbol{x}}(e^{\sum_{j=k+2}^n V(x_j)}|\boldsymbol{\mathcal{F}}_{k+1}))\\
         &\le \|f\|_\infty \|Q_{n-k-1}^V\mathbf{1}\|_\infty \mathbb{E}_{\boldsymbol{x}} (e^{\sum_{i=1}^k V(x_i)}\mathbb{E}_{\boldsymbol{x}}(\mathbf{1}_{C_k}e^{V(x_{k+1})}|\boldsymbol{\mathcal{F}}_k))
    \end{align*}
    In view of \eqref{conditionprob}, almost surely
    \[\mathbb{E}_{\boldsymbol{x}}(\mathbf{1}_{C_k}e^{V(x_{k+1})}|\boldsymbol{\mathcal{F}}_k)\le e^{\sup V} \mathbb{E}_{\boldsymbol{x}}(\mathbf{1}_{C_k}|\boldsymbol{\mathcal{F}}_k)\le C_V q^k d(x,x').\]
    (Throughout this proof the constant $C_V$ may change from line to line.) Therefore
    \begin{align*}
         \left|\mathbb{E}_{\boldsymbol{x}} (\mathbf{1}_{C_k}f(x_n)e^{\sum_{i=1}^n V(x_i)})\right|&\le C_V \|f\|_\infty \|Q_{n-k-1}^V \mathbf{1}\|_\infty q^k d(x,x') \mathbb{E}_{\boldsymbol{x}} e^{\sum_{i=1}^k V(x_i)}\\
         &\le C_V \|f\|_\infty \|Q_{n-k-1}^V \mathbf{1}\|_\infty \|Q_k^V \mathbf{1}\|_\infty q^k d(x,x'),
    \end{align*}
    The same estimate holds for the other term in $J_n^{\bad,k}$.
    Taking $\|Q_{n-k-1}^V \mathbf{1}\|_\infty\le \lambda_1^{n-k-1} M_n$ into account, we obtain (note that $\lambda_1$ is determined by $V$)
    \begin{equation}
        \sum_{k=0}^{n-1} |J_n^{\bad,k}|\le C_V\|f\|_\infty \sum_{k=0}^{n-1}\lambda_1^{n-k} \|Q_k^V \mathbf{1}\|_\infty q^k M_n d(x,x').\label{UF inequality2}
    \end{equation}

    \medskip

    Next we estimate $J_n^{\good}$. Noticing that $d(x_k,x'_k)\le q^k d(x,x')$ on $B_n$, almost surely
    \[\mathbf{1}_{B_n}\sum_{k=1}^n |V(x_k)-V(x_k')|\le C \Lip(V) \sum_{k=1}^\infty q^k d(x,x')\le C_Vd(x,x').\]
    This implies
    \begin{align*}
        \mathbb{E}_{\boldsymbol{x}} (\mathbf{1}_{B_n}|e^{\sum_{i=1}^n V(x_i)}-e^{\sum_{i=1}^n V(x_i')}|)&=\mathbb{E}_{\boldsymbol{x}}(\mathbf{1}_{B_n}|e^{\sum_{i=1}^n (V(x_i)-V(x_i'))}-1|e^{\sum_{i=1}^n V(x_i')})\\
        &\le C_V d(x,x') \mathbb{E}_{\boldsymbol{x}}e^{\sum_{i=1}^n V(x_i')}\le C_V \|Q_n^V \mathbf{1}\|_\infty d(x,x').
    \end{align*}
    It is thus straightforward to find
    \begin{align}
        |J_n^{\good}|&\le |\mathbb{E}_{\boldsymbol{x}} (\mathbf{1}_{B_n} (f(x_n)-f(x_n'))e^{\sum_{i=1}^n V(x_i)})|+|\mathbb{E}_{\boldsymbol{x}} (\mathbf{1}_{B_n} f(x_n')(e^{\sum_{i=1}^n V(x_i)}-e^{\sum_{i=1}^n V(x_i')}))|\notag\\
        &\le \Lip(f)\mathbb{E}_{\boldsymbol{x}} (\mathbf{1}_{B_n}d(x_n,x_n')e^{\sum_{i=1}^n V(x_i)})+\|f\|_\infty \mathbb{E}_{\boldsymbol{x}} (\mathbf{1}_{B_n} |e^{\sum_{i=1}^n V(x_i)}-e^{\sum_{i=1}^n V(x_i')}|)\notag\\
        &\le \Lip(f) q^n d(x,x') \mathbb{E}_{\boldsymbol{x}} e^{\sum_{i=1}^n V(x_i)}+C_V\|f\|_\infty \|Q_n^V \mathbf{1}\|_\infty d(x,x')\notag\\
        &\le C_V \|f\|_L\|Q_n^V \mathbf{1}\|_\infty d(x,x').\label{UF inequality3}
    \end{align}

    \medskip

    Gathering \eqref{UF inequality1}-\eqref{UF inequality3} together, we find that
    \[L_n\le C_V\|f\|_L \left(\sum_{k=0}^\infty \frac{\|Q_k^V \mathbf{1}\|_\infty}{\lambda_1^k} q^k+1\right)M_n.\]
    The infinite series is convergent owing to \eqref{geoserie}. Thus the proof is complete.
    \end{proof}

    It is worth mentioning that the last step of proof only requires the fixing rate $q$ to be strictly less than $1$, without extra smallness assumption on $q$.

    \medskip

    \textit{Step 2: Boundedness of $M_n$.} We start with finding an ``eigenmeasure" of $Q^{V*}$.

    \begin{lemma}\label{lem_lambda2}
        Assume the hypothesis $(\mathbf{I})$ is satisfied, and the function $V\in L(X)$ is given. Then there exists a constant $\lambda_2>0$ and $\mu\in \mathcal{P}(X)$ such that
        \begin{equation}
            Q^{V*} \mu=\lambda_2 \mu.\label{eigenmeasure}
        \end{equation}
        Moreover, $\supp(\mu)=X$ and $\lambda_2\le \lambda_1$ (recall $\lambda_1$ is defined by \eqref{lambda1definition}).
    \end{lemma}

    \begin{proof}[\rm \bf Proof of Lemma~\ref{lem_lambda2}]
        Consider the continuous map from $\mathcal{P}(X)$ into itself:
        \[\sigma\mapsto \frac{Q^{V*}\sigma}{(Q^{V*} \sigma)(X)}.\]
        Since $X$ is compact, so is $\mathcal{P}(X)$. Then Schauder's fixed point theorem implies that there exists a fixed point $\mu\in \mathcal{P}(X)$. In other words, let $\lambda_2=(Q^{V*}\mu)(X)>0$, we get \eqref{eigenmeasure}.
        
        Next we address $\supp(\mu)=X$, which involves hypothesis $(\mathbf{I})$. Equivalently speaking, for any $x\in X$ and $\varepsilon>0$, we need to show $\mu(B(x,\varepsilon))>0$. By hypothesis $(\mathbf{I})$, there exists $N\in \mathbb{N}$ and $p>0$ such that $P_N(y,B(x,\varepsilon))>p$ for any $y\in X$. Hence
        \begin{align*}
            \mu(B(x,\varepsilon))&=\lambda_2^{-N} (Q_N^{V*}\mu)(B(x,\varepsilon))=\lambda_2^{-N} \int_X Q_N^V(y,B(x,\varepsilon))\, \mu (dy)\\
            &\ge \lambda_2^{-N} \int_X e^{N\inf V}P_N(y,B(x,\varepsilon))\, \mu (dy)\ge \lambda_2^{-N}e^{N\inf V}p>0.
        \end{align*}

        Finally, noticing that
        \[\|Q_n^V \mathbf{1}\|_\infty\ge \langle Q_n^V \mathbf{1},\mu\rangle=\langle \mathbf{1},Q_n^{V*}\mu\rangle=\lambda_2^n\langle \mathbf{1},\mu\rangle=\lambda_2^n,\]
        we can deduce $\lambda_1\ge \lambda_2$ from its definition \eqref{lambda1definition}.
    \end{proof}

    Now we combine Lemma~\ref{lem_Ln<=Mn} and Lemma~\ref{lem_lambda2} to demonstrate the boundedness.
    \begin{lemma}\label{lem_Q_n^V1approxlambda^n}
        Assume the hypotheses $(\mathbf{I})$ and $(\mathbf{C})$ are satisfied, then for any $V\in L(X)$, the rescaled $L^\infty$ norm $M_n$ defined by \eqref{def_Mn} is bounded.
    \end{lemma}

    \begin{proof}[\rm \bf Proof of Lemma~\ref{lem_Q_n^V1approxlambda^n}]
        Assume this assertion fails, then there exists a subsequence $n_j\to \infty$ such that
        \[M_{n_j}=\lambda_1^{-n_j} \|Q_{n_j}^V \mathbf{1}\|_\infty\to \infty.\]
        According to Lemma~\ref{lem_Ln<=Mn} (with $f=\mathbf{1}$), we can rewrite $L_{n_j}\le C_V M_{n_j}$ as
        \[\Lip (f_j)\le C_V,\quad \text{where }f_j:=\|Q_{n_j}^V \mathbf{1}\|_\infty^{-1} Q_{n_j}^V \mathbf{1}.\]
        In particular, $f_j$ is uniformly bounded and equicontinuous. By Arzel\`a--Ascoli theorem, up to a subsequence (still labeled by $j$ for simplicity), $f_j$ converges uniformly to some $F\in C(X)$. Then
        \[\langle F,\mu\rangle=\lim_{j\to \infty} \|Q_{n_j}^V \mathbf{1}\|_\infty^{-1} \langle Q_{n_j}^V \mathbf{1},\mu\rangle=\lim_{j\to \infty} M_{n_j}^{-1} \lambda_1^{-n_j} \langle \mathbf{1},Q_{n_j}^{V*}\mu\rangle=\lim_{j\to \infty} M_{n_j}^{-1} (\lambda_2/\lambda_1)^{n_j}=0.\]
        Here we make use of $M_{n_j}\to \infty$ as well as $\lambda_1\ge \lambda_2>0$. Since $f_{n_j}\ge 0$ and $\|f_{n_j}\|_\infty=1$, we get $F\ge 0$ and $\|F\|_\infty=1$. This contradicts to $\langle F,\mu\rangle=0$ and $\supp(\mu)=X$.
    \end{proof}

    We can now conclude the proof of Proposition~\ref{prop_UF}.
    
    \begin{proof}[\rm \bf Proof of Proposition~\ref{prop_UF}]
        Gathering \eqref{Mn>=lambdan}, \eqref{Ln<=Mninequality}, and Lemma~\ref{lem_Q_n^V1approxlambda^n}, we obtain
        \[\Lip(Q_n^V f)=\lambda_1^n L_n\le C_V \|f\|_L \lambda_1^n M_n\le C_V \|f\|_L \lambda_1^n\le C_V \|f\|_L \|Q_n^V \mathbf
        {1}\|_{\infty}.\]
        This is exactly the desired inequality \eqref{UF inequality}.
    \end{proof}

    So far, the proof of Theorem~\ref{Thm LDP} is complete, through the chain of implications:
    \[\text{Prop~\ref{prop_UF}}\overset{\text{Thm~\ref{Thm Feynman-Kac}}}{\Longrightarrow}\text{Lem~\ref{lemma_FK}}\Longrightarrow\text{Lem~\ref{lemma_verifyKifer}}\overset{\text{Thm~\ref{Thm Kifer}}}{\Longrightarrow}\text{Thm~\ref{Thm LDP}}.\]

    \section{Further discussion}\label{sec_discussion}

    Although we have successfully established the Donsker--Varadhan LDP for various randomly forced PDEs, it remains unknown whether LDP still holds for initial data $u_0\not \in \supp(\mu_*)$. To the best of our knowledge, this question is open even for parabolic PDEs with nondegenerate noise; see, e.g., \cite{JNPS-15,JNPS-18,JNPS-21}. The method of our work relies heavily on the uniform irreducibility, which guarantees the rate function $I$ to be independent of the initial state. And uniform irreducibility unavoidably fails without restricting to $\supp(\mu_*)$: due to the dissipation of energy and boundedness of noise, the energy of solution has an apriori bound.

    Even for very simple finite-state Markov chains, the lack of irreducibility may destroy both the convexity of rate function, and its independence on the initial state; see, e.g., \cite[remark after Theorem 3.1.2]{DZ10}. If the rate function is non-convex, then we cannot characterize it as the Legendre transform of pressure function. And if it depends on the initial data, then the interaction among different initial states can be complicated.

    Very recently, we realized that when the noise is suitably unbounded (e.g., white noise), then $\supp(\mu_*)$ coincides with the whole phase space, and thus $u_0\in \supp(\mu_*)$ is no longer restrictive. Therefore, we believe our methods in \cite{CLXZ24} and this paper to be applicable to establish the LDP for stochastic nonlinear wave equations (i.e.~\eqref{waveequation} with $\eta$ standing for white-in-time noise) with respect to any initial data in $\mathcal{H}$, extending the local LDP result in \cite{MN-18}. The techniques are more intricate due to the white noise. We plan to address this question in a subsequent work.

    \vspace{5mm}
    \noindent\textbf{Acknowledgments\addcontentsline{toc}{section}{Acknowledgments}} \; The authors thank the anonymous referee for insightful comments and suggestions. We are also grateful to Ziyu Liu for valuable discussions and advices during the preparation of the paper. Yuxuan Chen would like to appreciate his PhD advisor Zhifei Zhang. Shengquan Xiang is partially supported by NSFC 12301562.

\begin{appendix}
    
    \section*{Appendix}
    
    \stepcounter{section}
    
    \numberwithin{equation}{subsection}
    \numberwithin{theorem}{subsection}
    
    \setcounter{equation}{0}
    \setcounter{theorem}{0}
    
    \renewcommand\thesubsection{\Alph{subsection}}

    \stepcounter{subsection}

    \hypertarget{Appendix_A}{}
    \noindent\Alph{subsection}. \,\textbf{Kifer's criterion for uniform LDP}
    
    Suppose $X$ is a compact metric space, $(\Theta,\prec)$ is a directed set, and $L_\theta\in \mathcal{P}(X)\, (\theta\in \Theta)$. Let the scale $r_\theta\in (0,\infty)$ satisfy $\lim_{\theta\in \Theta} r_\theta=\infty$. Assume the convergence of
    \[\Lambda(V):=\lim_{\theta\in \Theta} \frac{1}{r_\theta} \log \mathbb{E} e^{r_\theta \langle V,L_\theta\rangle}\quad \text{for }V\in C(X).\]
    The rate function $I\colon \mathcal{P}(X)\to [0,\infty]$ is defined by Legendre transform \eqref{legendre}. And conversely \eqref{dualrelation} is valid. The equilibrium state is defined as in Definition~\ref{def_es}.

    \begin{theorem}[{\cite[Theorem 2.1]{Kifer-90}}]\label{Thm Kifer}
        Assume $\Lambda\colon C(X)\to \mathbb{R}$ is well-defined, and that there exists a dense linear subspace $\mathcal{V}\subset C(X)$, such that the equilibrium state is unique for any $V\in \mathcal{V}$. Then the LDP holds for $L_\theta$ with respect to the convex good rate function $I$.
    \end{theorem}

    \hypertarget{Appendix_B}{}
    \stepcounter{subsection}
    \noindent\Alph{subsection}. \,\textbf{Asymptotics of generalized Markov semigroups}

    Let $X$ be a compact metric space. A generalized Markov kernel consists of a family $\{Q(x,\cdot):x\in X\}\subset \mathcal{M}_+(X)$ such that $0<Q(x,X)<\infty$ for any $x\in X$, and the map $x\mapsto Q(x,\cdot)$ is continuous from $X$ to $\mathcal{M}_+(X)$. We denote by $Q_n(x,\cdot)$ the $n$-fold iteration of $Q(x,\cdot)$. With a slight abuse of notation, let the operators
    \begin{equation*}
        Q_n\colon C(X)\rightarrow C(X),\quad Q_n^*\colon \mathcal{M}_+(X)\rightarrow \mathcal{M}_+(X)
    \end{equation*}
    be defined by
    \begin{equation*}
    Q_nf(x)=\int_X f(y)\, Q_n(x,dy),\quad Q_n^*\sigma(A)=\int_X Q_n(x,A)\, \sigma(dx),
    \end{equation*}
    where $f\in C(X), \sigma\in\mathcal{M}_+(X)$ and $A\in\mathcal{B}(X)$. It can be verified that $Q_n$ and $Q_n^*$ are semigroups.
    
    \begin{definition}[Determining family]\label{def_determine}
        A family $\mathcal{C}\subset C(X)$ is called determining if the following property holds: given any $\mu,\nu\in\mathcal{M}_+(X)$, if $\langle f,\mu\rangle=\langle f,\nu\rangle$ for any $f\in \mathcal{C}$, then $\mu=\nu$.
    \end{definition}
        
    \begin{theorem}[{\cite[Theorem 2.1]{JNPS-15}}]\label{Thm Feynman-Kac}
        Under the above settings, assume further that:
        \begin{itemize}
            \item[\tiny$\bullet$] {\bf Uniform irreducibility.} For any $\varepsilon>0$, there exists $N\in \mathbb{N}$ and $p>0$, such that
            \[Q_N(x,B(y,\varepsilon))>p,\quad \text{for any } x,y\in X.\]
               
            \item[\tiny$\bullet$] {\bf Uniform Feller property.} There exists a determining family $\mathcal{C}\subset C(X)$, such that for any $f\in \mathcal{C}$, the sequence $\{\|Q_n f\|^{-1}_\infty Q_nf:n\in \mathbb{N}\}$ is equicontinuous.
        \end{itemize}
        Then there exists $\lambda>0$, $h\in C(X)$ and $\mu\in\mathcal{P}(X)$, such that for any $f\in C(X)$ and $\nu\in\mathcal{M}_+(X)$,
        \begin{align*}
            &Qh=\lambda h,\quad Q^*\mu=\lambda\mu,\quad \langle h,\mu\rangle =1,\quad h>0,\quad \supp(\mu)=X,\\
            &\lambda^{-n}Q_nf\rightarrow\langle f,\mu\rangle h\quad\text{in }C(X),\quad \lambda^{-n}Q_n^*\nu\rightarrow\langle h,\nu\rangle \mu\quad \text{in }\mathcal{M}_+(X).
        \end{align*}
    \end{theorem}

\end{appendix}

    \normalem
    \bibliographystyle{plain}
    \bibliography{References}
\end{document}